\documentclass[12pt,twoside]{amsart}

\usepackage{enumerate,latexsym,amsthm}

\usepackage{amsmath}
\usepackage{amssymb}

\theoremstyle{plain}

\newtheorem{proposition}{Proposition}[section]
\newtheorem{theorem}[proposition]{Theorem}
\newtheorem{lemma}[proposition]{Lemma}
\newtheorem{corollary}[proposition]{Corollary}

\theoremstyle{definition}
\newtheorem{definition}[proposition]{Definition}

\newtheorem{remark}[proposition]{Remark}

\def\A{\mathcal{A}}
\def\C{\mathbb{C}}
\def\F{\mathcal{F}}

\def\L{\mathcal{L}}

\def\O{\mathcal{O}}
\def\P{\mathbb{P}}

\def\Z{\mathbb{Z}}

\title{
Syzygies, regularity and toric varieties
} 
\author{Milena Hering}
\date{\today}
\begin{document} 
\begin{abstract}
Let $\A$ be an ample line bundle on a projective toric
variety $X$ of dimension $n$. 
We show that if  $\ell \geq n-1+p$,  then $\A ^{\ell}$  satisfies
the property $N_p$. 
Applying similar methods, we obtain a combinatorial theorem:
For a given lattice polytope $P$ we give a criterion
for an integer $m$ to guarantee that $mP$ is normal.\end{abstract}

\maketitle 
 \section{Introduction}
Let $\L$ be a globally generated, ample line bundle on a projective variety $X$ 
over a field $k$ of characteristic zero 
and 
let $\phi _{\L}$ be the map of $X$ into $\P (H^0(X,\L))$. 
Let $S$ be the symmetric algebra $\mathrm{Sym}^{\bullet}H^0(X,\L)$
and let $R=\bigoplus_{m}H^0(X,\L ^m)$, a finitely
generated S-module. 

Let 
$$ 0\to E_k \to \cdots \to E_1 \to E_0 \to R \to
0$$
be a minimal free graded resolution of the $S$-module $R$. 
$E_i$ is also called the {\it $i$'th syzygy module} of $R$.

We say that an ample line bundle 
$\L$ {\it satisfies property} $N_0$ if $E_0=S$, 
and it {\it satisfies property} $N_p$ if 
$E_0=S$ and $E_i \cong \bigoplus S(-i-1)$ for $1\leq i\leq p$.

For example, if $\L$ satisfies $N_0$, $\L$ is very 
ample and if $X$ is normal, $\phi_{\L}$ embeds 
$X$ into $\P (H^0(X,\L)$ as a projectively normal variety. If it
also satisfies $N_1$, it has quadratic defining 
equations and $N_2$ implies that the relations 
among these equations are linear.  
For $p=0, 1$, the property $N_p$ has been studied by 
Mumford in 
\cite{M}, who called such line bundles {\it normally generated}
and {\it normally presented}, respectively. For a survey of this
property, we refer to \cite{PAG}, Section 1.8.D.

Recall that
for an
ample and globally generated 
line bundle $\A$ on $X$, 
a sheaf $\F$ is called {\it $m$-regular with respect to $\A$} in 
the sense of Castelnuovo-Mumford, if 
\begin{equation*} H^i(X,\F \otimes \A^{m-i}) = 0 
\textrm{ for all }i\geq 1.
\end{equation*}
\begin{definition}
We will call $\A$ \emph{$m$-autoregular}, if it is $m$-regular
with respect to itself, i.e., if
\begin{equation*}
H^i(X,\A^{m+1-i}) = 0
\textrm{ for all }i\geq 1.
\end{equation*}
\end{definition}

The following theorem is due to Gallego and Purnaprajna. 
We include a short proof in this paper. 

\begin{theorem}[\cite{GP}, Theorem 1.3.]\label{thetheorem}
Let $\A$ be an ample line bundle on X that is globally generated, 
and suppose that $\A$ is $m$-autoregular for some integer $m$. 
If $\ell\geq \max\{1,m+p\}$ and $p\geq 1$, 
then $\A ^{\ell}$ satisfies property $N_p$. 
\end{theorem}

This theorem has a particularly nice application to an ample line
bundle $\A$ on a  projective 
toric varietiy $X$ of dimension $n\geq 2$. 
Every ample line bundle $\A$ on a toric variety 
of dimension $n$
is $(n-1)$-autoregular, since it is 
globally generated and since the higher 
cohomology of a globally generated
line bundle 
vanishes.
 Hence we can apply Theorem \ref{thetheorem} to $\A$. Together
 with the fact that if $\ell \geq n-1$,  then
 $\A ^{\ell}$ satisfies property $N_0$ (\cite{EW}, 
 \cite{LTZ}, \cite{NO}),
we 
obtain the following Corollary.
\begin{corollary}\label{th:1}
Let $\A$ be an ample line bundle on a toric variety $X$ of
dimension $n\geq 2$. Then
$\A ^{\ell}$ satisfies property $N_p$ when $\ell\geq n-1+p$
and $p\geq 0$. 
\end{corollary}

This result has first appeared in a preprint by Hal Schenck
and Gregory Smith \cite{S}. We can fix a gap in their
proof by applying Theorem \ref{thetheorem}.

The proofs of the case $p=0$ in 
\cite{EW}, \cite{LTZ} and \cite{NO} use the combinatorics
of lattice points in polytopes.
Nakagawa and Ogata prove 
the case $p=1$ following Mumford in \cite{M}; 
in \cite{BGT}, 
Bruns, Gubeladze and Trung prove the case $p=0$ and $p=1$
in a special case using the commutative algebra associated
to polytopal semigroup rings.

Koelman gives a combinatorial 
criterion for an ample line bundle on a toric 
surface to be normally
presented in \cite{K2}:
Let $\L = \O _X(D)$ for some torus invariant Cartier divisor 
$D$ and let $P_D$ be the polytope associated to $D$. Recall that
$P_{mD} = mP_D$.
Then Koelman proves that 
$\O _X(D)$  is normally presented if $P_D$ contains 
more than 3 lattice points in its boundary; in particular
he shows that 
for any ample line bundle $\A$ on a toric surface, $\A ^2$
is normally presented. 

Investigating the regularity of ample line bundles on 
toric varieties, we obtain an interesting
corollary for lattice polytopes.  
Recall that a lattice polytope is \emph{normal} 
if every lattice point in $mP$ is the sum of $m$ lattice
points in $P$. Note that on a 
toric variety $X$, $\O _X(D)$ is normally generated
if and only if $P_D$ is normal. 
Moreover, let 
$d(P)$ be the largest integer such that 
$d(P) P$ does not contain any lattice points in its relative 
interior. 

\begin{corollary}\label{th:3}
Let $V$ be a real vector space of dimension $n$, and let
$M \subset V$ be a lattice of full rank. Let $P$ 
be a lattice polytope of dimension $n$. 
Then $\ell P$ is a normal polytope for all 
$\ell\geq \max\{n-d(P),1\}$.
\end{corollary}

In section \ref{se:preliminaries} 
we recall the cohomological criterion of Green, 
Ein and Lazarsfeld 
(\cite{EL}) for 
a line bundle to satisfy the property $N_p$, and 
some basic facts about the regularity
of coherent sheaves. 
In  section \ref{se:proof} we will use these facts
to prove Theorem \ref{thetheorem}. In the last section we 
give a combinatorial condition for  
the regularity of ample line bundles on toric varieties, 
and in this way we obtain some better bounds than in 
Corollary \ref{th:1}.

I have learnt most of the 
techniques and theorems 
from the preprint 
by Hal Schenck and Gregory Smith \cite{S};
in particular, 
the idea to use the regularity of powers
of the 
vector bundles $M _{\L}$ defined in (\ref{M}), which
had also been applied by Gallego and Purnaprajna in \cite{GP}.
Similar methods appear already in \cite{EL}. 

I wish to thank W. Fulton and R. Lazarsfeld  for 
helpful discussions related  to this work, 
and also A. Bayer 
for comments on earlier versions of this note.

\section{Preliminaries}\label{se:preliminaries}
To a globally generated line bundle $\L$ we associate 
a vector bundle $M_{\L}$, defined by the following exact sequence:
\begin{equation}\label{M} 
0\to M_{\L} \to H^0(X,\L)\otimes O_X \to \L \to 0.
\end{equation}
The following cohomological criterion for the 
property 
$N_p$ is the main tool in the proof of 
the theorem. It can be found in \cite{EL}. 
\begin{theorem}\label{cor:1}
Let $\mathrm{char}(k)=0$ and $\L$ be an ample line bundle
that is globally generated.
Then $\L$ satisfies property $N_p$ if
\begin{equation*}
H^1(X,M_{\L}^{\otimes k}\otimes \L ^j)= 0 \ \ \mathrm{ for }\ \ 
0\leq k \leq p+1\ \  \mathrm{ and } \ \ j\geq 1.
\end{equation*}
\end{theorem}

Moreover, we will need the following properties of regular sheaves.
\begin{remark}\label{rem:1}
It follows from the definition that when $\F$ is $m$-regular
with respect to $\A$, then $\F \otimes \A ^d$ is $(m-d)$-regular
with respect to $\A$ for all $d\in \Z$.
\end{remark}
\begin{lemma}\label{sesregularity}
Let $\A$ be an ample and globally generated
line bundle on $X$ 
and let 
\begin{equation}\label{eq:ses}
0\to \F ' \to \F \to \F '' \to 0
\end{equation}
be a short exact sequences of coherent sheaves on $X$. 
Suppose that 
$\F$ is $r$-regular with respect to $\A$, that  
$\F ''$ is $(r-1)$-regular with respect to
$\A$, 
and that 
the map of global sections
\begin{equation}\label{eq:3}
H^0(X,\F \otimes \A ^{r-1}) \to H^0(X,\F '' \otimes
\A ^{r-1})
\end{equation}
is surjective.
Then $\F '$ is $r$-regular with respect to $\A$.
\end{lemma}
This follows from twisting 
(\ref{eq:ses}) with appropriate powers of $\A$ and 
applying the definition of regularity to the long exact sequence
of cohomology groups. 

\begin{theorem}[Mumford's Theorem, 
\cite{PAG} Theorem 1.8.5.]
\label{MT}
Let $\A$ be an ample and globally generated
line bundle on a projective variety
$X$. 
Let $\F$ be a sheaf on $X$ that is $m$-regular with respect
to $\A$. Then for every $\ell\geq 0$:

\begin{enumerate}
\item\label{M:1} $\F\otimes\A^{m+\ell}$ is generated by its
global sections.
\item\label{M:2}The natural maps 
\begin{equation*}
H^0(X,\F\otimes \A^{m})\otimes H^0(X,\A^{\ell})
\to H^0(X,\F\otimes \A^{m+\ell})
\end{equation*}
are surjective.
\item\label{M:3} $\F$ is $(m+\ell)$-regular with respect to $\A$.
\end{enumerate}
\end{theorem}

\section{
The relation between regularity and the property $N_p$}
\label{se:proof}
\begin{proof}[Proof of Theorem \ref{thetheorem}]
Let $\L = \A ^{\ell}$ with $\ell\geq \max \{1,m+1\}$. 
Since $\L$ is ample and globally generated,
we can associate a vector bundle $M_{\L}$ to it as in 
(\ref{M}).

We 
claim that $M_{\L}^{\otimes k}$ is 
$(m+k)$-regular with respect to $\A$ for $k\geq 1$.
Granting this claim, we see that in particular,
$H^1(X, M_{\L}^{\otimes k}\otimes \A ^d) = 0$ for $d\geq
m+k-1$. Therefore, when $\ell\geq \max\{m+p,1\}$ and $j\geq 1$,
$j\ell \geq \ell\geq m+p \geq m+k-1$ for $p+1\geq k\geq 0$
and hence $H^1(X, M_{\L}^{\otimes k}\otimes \L^{j}) = 0$ in
this case.
When $k=0$, the vanishing  follows from the
$m$-regularity of $\A$. It follows from Theorem \ref{cor:1} 
that $\L$ satisfies the 
property $N_p$.

We will prove the claim by induction on $k$,  
applying  Lemma \ref{sesregularity} with $r=m+k$ to the short exact sequence
(\ref{M}) defining $M_{\L}$ twisted by $M_{\L}^{\otimes (k-1)}$.

The case $k=1$ is a special case of Lemma 3.1. of
\cite{AK}, but we can see it directly as
follows. Observe that since $\ell -1\geq m$,  
$\A$ is $(\ell -1)$-autoregular and so  
we can apply Mumford's theorem
 (\ref{M:2}) to $\A $ 
to see 
that for $\L = \A \otimes \A ^{\ell-1}$, 
the map of global sections
\begin{equation*}
H^0(X,\L)\otimes H^0(X,\A^{m}) \to H^0(X,\L \otimes 
\A^{m})
\end{equation*}
is surjective.
Using Remark \ref{rem:1} and the fact that cohomology commutes 
with tensoring with a vector space, 
we see that $H^0(X,\L)\otimes \O _X$ is $(m+1)$-regular 
with respect to $\A$. Similarly,  
$\L$ is $(m-\ell+1)$-regular, hence $m$-regular with respect to
$\A$ 
since $\ell\geq 1$.
Now Lemma \ref{sesregularity} implies that $M_{\L}$ is $(m+1)$-regular with
respect to $\A$.

For $k> 1$, we substitute $\A^{\ell}$ 
for $\L$ in (\ref{eq:3}) and apply (\ref{M:2}) in
Mumford's theorem to $M_{\L}^{\otimes (k-1)}$, which
is $(m+k-1)$-regular by the induction hypothesis, to see that
the map
of global sections
\begin{equation*}
\begin{split}
\shoveleft{
H^0(X,M_{\L}^{\otimes (k-1)}\otimes \A^{m+k-1})\otimes 
H^0(X,\A^{\ell})  }
\\ 
\shoveright{ \to H^0(X,
M_{\L}^{\otimes (k-1)}\otimes \A^{m+k-1+\ell})}
\end{split}
\end{equation*}
is surjective. 
Moreover, by the induction hypothesis $M_{\L}^{\otimes (k-1)}
\otimes H^0(X,\L)$ is $(m+k-1)$-regular, in particular
it is $(m+k)$-regular; similarly 
$M_{\L}^{\otimes (k-1)}\otimes \L$
is $(m+k-1)$-regular, since $\ell \geq 1$. 
Now the claim follows from 
Lemma \ref{sesregularity}. \end{proof}

\begin{remark}
That $\A ^{m+1}$ satisfies property $N_0$ follows directly
from Mumford's theorem (\ref{M:2}).
\end{remark}

\section{The regularity of ample line bundles on toric varieties
}
In this section we will use a theorem by David Cox and Alicia
Dickenstein to compute the regularity of a given line bundle
on a toric variety in terms of combinatorial properties of 
a polytope associated to the line bundle.

\begin{theorem}[\cite{CD}, Theorem 1.3.]\label{th:CD}
Let $D$ be a torus invariant Cartier divisor on a complete
toric variety $X$ and let $P_D$ be the polytope associated
to $D$. Assume that $\O _X(D)$ is globally generated. 
Then 
\begin{enumerate}
\item 
$H^i(X, \O _X(-D)) = 0$ for all $i\neq \dim (P_D)$.
\item
There is an isomorphism
\begin{equation*}
H^{\dim (P_D)}(X,\O _X(-D))\cong \bigoplus _{u\in \mathrm{relint}
(P_D)\cap M}\C\chi ^{-u}
\end{equation*}
that is equivariant with respect to the torus action.
\end{enumerate}
\end{theorem}
In particular, if $P_D$ has no interior lattice points, the
cohomology of $\O _X(-D)$ vanishes. This motivates
the definition 
of $d(P)$ for a polytope $P$ in the introduction.

\begin{definition}
Let $\A$ be an ample line bundle that is globally
generated. We define the $autoregularity$ of $\A$ to 
be the smallest intger $m$ such that $\A$ is 
$m$-autoregular.
\end{definition}
\begin{proposition}
Let $D$ be a torus invariant ample divisor on a complete 
toric variety $X$ of dimension $n$
and let $P_D$ be the polytope associated to 
$D$. Let $m= n-1-d(P_D)$. 
Then the autoregularity of $\O _X(D)$ is $m$.
\end{proposition}

\begin{proof}
We have to show that, for $i\geq 1$,
\begin{equation*}
H^i(X,\O _X( (m+1-i)D)) = 0.
\end{equation*}
First recall that on a toric variety for any ample divisor
$D$, $\O _X(D)$ is globally generated. So 
when $m+1-i\geq 0$, the statement follows from the vanishing
of the higher cohomology of globally generated line bundles
on any complete toric variety. 
(See for example \cite{F}, section 3.4.).

When $m+1-i < 0$, $-(m+1-i)D$ is a positive integer 
multiple of $D$, 
in particular it is 
globally
generated and we can apply Theorem \ref{th:CD}. Since 
the dimension of the polytope associated to an ample divisor 
equals the dimension of the variety, the assertion follows
for 
$i\neq n$.
When $i=n$, $m+1-n  
=-d(P_D) $.
But $d(P_D) P_D$ does not contain any interior 
lattice points by definition, 
and so Theorem \ref{th:CD} (2) implies that also 
$H^n(X, \O _X(-d(P_D) D) = 0$.
Moreover $H^n(X, \O _X(-(d(P_D) +1)D))\neq 0$, 
since $( d(P_D)+1)P$ does contain interior lattice 
points,
so $\O _X (D)$ is not $(m-1)$-autoregular.
\end{proof}

Applying Theorem \ref{thetheorem} to 
an ample divisor $D$ on a toric variety $X$, we obtain the 
following Corollary, which implies
Corollary \ref{th:3}.
\begin{corollary}\label{co:treg}
Let $D$ be a torus invariant ample divisor on a complete
toric variety $X$ of dimension $n$,
let $P_D$ be the polytope associated to
$D$ and let $m= n-1-d(P_D)$.
If $p\geq 1$ and $\ell\geq \max \{m+p,1\}$, 
then 
the line bundle $\O _X(\ell D)$ satisfies
property $N_p$. 
In particular $(m+1)D$ is very ample and normally generated. 
\end{corollary}

\bibliographystyle{alphaspecial}
\bibliography{all}

\begin{thebibliography}{Mum70}

\bibitem[AK02]{AK}
Donu Arapura and Dennis Keeler.
\newblock Frobenius amplitude and strong vanishing theorems for vector bundle.
\newblock 2002.
\newblock math.AG/0202129.

\bibitem[BGT97]{BGT}
Winfried Bruns, Joseph Gubeladze, and Ng{\^o}~Vi{\^e}t Trung.
\newblock Normal polytopes, triangulations, and {K}oszul algebras.
\newblock {\em J. Reine Angew. Math.}, 485:123--160, 1997.

\bibitem[CD03]{CD}
David Cox and Alicia Dickenstein.
\newblock Vanishing and codimension theorems for complete toric varieties.
\newblock 2003.
\newblock math.AG/0310108v1.

\bibitem[EL93]{EL}
Lawrence Ein and Robert Lazarsfeld.
\newblock Syzygies and {K}oszul cohomology of smooth projective varieties of
  arbitrary dimension.
\newblock {\em Invent. Math.}, 111(1):51--67, 1993.

\bibitem[EW91]{EW}
G{\"u}nter Ewald and Uwe Wessels.
\newblock On the ampleness of invertible sheaves in complete projective toric
  varieties.
\newblock {\em Results Math.}, 19(3-4):275--278, 1991.

\bibitem[Ful93]{F}
William Fulton.
\newblock {\em Introduction to toric varieties}, volume 131 of {\em Annals of
  Mathematics Studies}.
\newblock Princeton University Press, Princeton, NJ, 1993.

\bibitem[GP99]{GP}
F.J. Gallego and B.P. Purnaprajna.
\newblock Projective normality and syzygies of algebraic surfaces.
\newblock {\em J. Reine Angew. Math.}, 506:145--180, 1999.

\bibitem[Koe93]{K2}
Robert~Jan Koelman.
\newblock A criterion for the ideal of a projectively embedded toric surface to
  be generated by quadrics.
\newblock {\em Beitr\"age Algebra Geom.}, 34(1):57--62, 1993.

\bibitem[Laz]{PAG}
Robert Lazarsfeld.
\newblock {\em Positivity in Algebraic Geometry}.
\newblock to appear.

\bibitem[LTZ93]{LTZ}
Ji~Yong Liu, Leslie~E. Trotter, Jr., and G\"unter~M. Ziegler.
\newblock On the height of the minimal hilbert basis.
\newblock {\em Results in Mathematics}, 23(3-4):374--376, 1993.

\bibitem[Mum70]{M}
David Mumford.
\newblock Varieties defined by quadratic equations.
\newblock In {\em Questions on Algebraic Varieties (C.I.M.E., III Ciclo,
  Varenna, 1969)}, pages 29--100. Edizioni Cremonese, Rome, 1970.

\bibitem[ON02]{NO}
Shoetsu Ogata and Katsuyoshi Nakagawa.
\newblock On generators of ideals defining projective toric varieties.
\newblock {\em Manuscripta Math.}, 108(1):33--42, 2002.

\bibitem[SS03]{S}
Hal Schenck and Gregory~G. Smith.
\newblock Syzygies of projective toric varieties.
\newblock 2003.
\newblock math.AG/0308205.

\end{thebibliography}

\end{document}